\documentclass[11pt]{amsart}
\usepackage{amssymb, amscd}
\usepackage{latexsym}
\usepackage[all,dvips]{xy}
\numberwithin{equation}{section}

\newcommand{\QQ}{\mathbb Q}
\newcommand{\CC}{\mathbb C}
\newcommand{\EE}{\mathbb E}
\newcommand{\FF}{\mathbb F}
\newcommand{\GG}{\mathbb G}
\newcommand{\ZZ}{\mathbb Z}
\newcommand{\PP}{\mathbb P}
\newcommand{\RR}{\mathbb R}
\newcommand{\VV}{\mathbb V}
\newcommand{\WW}{\mathbb W}

\newtheorem{theorem}{Theorem}[section]

\newtheorem{conjecture}[theorem]{Conjecture}
\newtheorem{definition-lemma}[theorem]{Definition-Lemma}

\theoremstyle{definition}

\theoremstyle{remark}

\begin{document}

\title[Modular Forms and the Cohomology of Local systems]{Exploring Modular Forms 
and the Cohomology of Local systems on Moduli Spaces by Counting Points}
\author{Gerard van der Geer}
\address{Korteweg-de Vries Instituut, Universiteit van
Amsterdam,
\newline Postbus 94248, 1090 GE Amsterdam, The Netherlands.}
\email{geer@science.uva.nl}

\begin{abstract} We report on a joint project in experimental 
mathematics with Jonas Bergstr\"om and Carel Faber 
where we obtain information about modular forms by 
counting curves over finite fields. 
\end{abstract}
\maketitle
\section{Introduction} \label{sec-intro}
According to Arnold mathematics is the part of physics where experiments
are cheap. I take it that he meant `cheap' in
the literal sense, that is, not requiring large investments of goods and money.
But experiments in mathematics may require a considerable 
investment of time, mental energy  and
computational power. On the other hand they can be extremely rewarding. 
Here we report on a  project of experimental mathematics concerning
moduli spaces and modular forms that has been going on for quite some time.
We used counts of points over finite fields to explore the cohomology
of local systems on moduli spaces of low genus curves (\cite{FvdG, BFG, BFGg3, JBvdG2}). 
As a side result
it produces a lot of information about Siegel modular forms of low degree.
This comes about since on the one hand
modular forms admit a cohomological description
as part of the cohomology of local systems on moduli spaces of abelian 
varieties and on the other hand
abelian varieties of low dimension allow a description 
in terms of curves. By using comparison theorems for the
cohomology we can gauge this cohomology by looking
at the situation over a finite field and listing all the isomorphism classes
of curves of this genus over the finite field and the order of their automorphism groups,
and by calculating the action of
Frobenius on the cohomology of the corresponding curve.

Birch made such counts for elliptic curves in the 1960s in  \cite{Birch} and from his formulas
one sees that one gets information about the traces of the Hecke operators on the spaces of cusp forms. 
In a joint project with
Carel Faber we started counting points of the moduli spaces 
${\mathcal M}_{2,n}$
of $n$-pointed curves of genus $2$ over finite fields
in order to probe the cohomology of these moduli spaces. 
We realized that knowledge of the cohomology of the moduli spaces ${\mathcal M}_{2,n}$
is equivalent to information about the cohomology of local systems on ${\mathcal M}_2$
and also realized that this way we could get information about vector-valued Siegel modular forms.  
It turned out to be wonderfully effective and our heuristic results have since then been
turned into theorems. In joint work with Bergstr\"om we extended
this to local systems on the moduli space ${\mathcal M}_2[2]$ of curves of genus $2$ 
with a level $2$ structure
and then later also to the case of ${\mathcal M}_3$, the moduli 
of curves of genus $3$,
where it makes predictions about Siegel modular forms of degree $3$. 
It was a surprise that this was feasible
and it thus created a new window on Siegel modular forms.
Apart from these cases we also looked at the case of Picard modular forms.
Another chapter is the case of modular forms on ${\mathcal M}_3$, 
called Teichm\"uller modular forms. We will not treat this last aspect
here  and have to 
refer to future papers.

One of the nice aspects of the theory of elliptic modular forms is the availability of concrete examples
and accessibility of the topic for doing experiments. Concrete examples for higher degree
Siegel modular forms were much harder to get by, especially for vector-valued modular
forms. But recently things have changed quite a lot. For example the work of Chenevier-Renard
and Ta\"{\i}bi (see \cite{C-R, Taibi}) where applying Arthur's results on the 
trace formula (even if all has not yet been proved) provided a wealth of new data on Siegel modular forms.
No doubt we can await many new surprises and discoveries in this beautiful corner of 
mathematical nature.

Finally I would like to thank my collaborators Jonas Bergstr\"om and Carel Faber 
for their comments and the continued pleasant cooperation in this project.
I also thank Fabien Cl\'ery with whom I worked on the purely 
modular forms side of the project.
\section{Modular Forms}\label{modforms}
Every mathematician knows or should know the definition of (elliptic) modular form: 
a holomorphic function $f: \mathfrak{H} \to {\CC}$ on the upper
half plane $\mathfrak{H}$  of ${\CC}$ satisfying (for a fixed $k\in {\ZZ}$)
$$
f((a\tau+b)(c\tau+d)^{-1})=(c\tau+d)^kf(\tau)
$$
for all $\tau \in \mathfrak{H}$ and $(\begin{smallmatrix} a & b \\ c & d \\ \end{smallmatrix}) 
\in {\rm SL}(2,{\ZZ})$,
in particular $f(\tau+1)=f(\tau)$, 
and thus admitting a Fourier series
$$
f=\sum_{n\in {\ZZ}} a(n)\,  q^n \qquad \hbox{\rm with $q=e^{2\pi i \tau}$}\, ,
$$
where we require additionally that $a(n)=0$ for $n<0$. The exponent $k\in {\ZZ}$ 
is called the weight of the modular form. 
Examples of modular forms can be found everywhere in mathematics, 
and in fact they appeared early on in our history, for example in the papers of Jacobi, 
but the notion of a modular form was formalized only
at the end of the nineteenth century, apparently by Klein who introduced
the word `Modulform.' 

Modular forms 
of given weight $k$ form a finite-dimensional complex vector space denoted by $M_k$ (or
$M_k({\rm SL}(2,{\ZZ}))$) and since the product of two modular forms of weight $k$
and $l$ is a modular form of weight $k+l$ they form a graded ring
and ${\CC}$-algebra
$$
M=\oplus_{k} M_k \, .
$$
A modular form vanishing at infinity, that is, with $a(0)=0$, is called a cusp form.
The subspace of cusp forms is denoted by $S_k$.
The cusp form of smallest non-zero weight is $\Delta=\sum_{n=1}^{\infty} \tau(n) \, q^n$ 
which can be defined by an infinite product
$$
q \prod_{n=1}^{\infty} (1-q^n)^{24}= 
q-24 \, q^2+252\, q^3-1472 \, q^4+4830 \, q^5 +\cdots
$$
and this turns out to be 
a cusp form of weight $12$ generating the ideal of cusp forms in the algebra $M$.

The structure of this ${\CC}$-algebra $M$ is known: it is a polynomial algebra
$$
M= {\CC}[E_4,E_6],
$$
where the modular forms $E_4$ and $E_6$ of weight $4$ and $6$ are Eisenstein series
defined for even $k\geq 4$ via
$$
E_{k}= \frac{1}{2} \sum_{(c,d)=1} (c\tau+d)^{-k}
$$
where $(c,d)$ runs over coprime pairs of integers and we assume $k\geq 4$ to get convergence; 
it has the Fourier series 
$$ 
1-\frac{2 k}{B_k}\sum_{n=1}^{\infty} \sigma_{k-1}(n)\, q^n
$$
where $\sigma_r(n)=\sum_{1\leq d|n} d^r$ and $B_k$ is the $k$th Bernoulli number. The cusp form
$\Delta$ equals the expression $(E_4^3-E_6^2)/1728$.

We can replace ${\rm SL}(2,{\ZZ})$ by a subgroup of finite index, like for example
the subgroup $\Gamma_0(N)$ where we require for a matrix $(a,b;c,d)$ that
$c$ is divisible by $N$. 
Then we find more modular forms. However, for the most part of this report
 we restrict to the case of level $N=1$.

Modular forms show up unexpectedly in many parts of mathematics, for example in the
theory of lattices as generating series
for the number of vectors of given norm in an even definite 
unimodular lattice. 
They show up ubiquitously 
if one studies the moduli of elliptic curves or other algebro-geometric
objects. 
But they also show up in mathematical physics.
As a rule the Fourier coefficients of these 
modular forms show interesting arithmetical properties and seem to contain 
magical information.  

It was Hecke who saw around 1937 how to extract this information 
from the Fourier series. He introduced the 
operators $T(n)$ for $n\in {\ZZ}_{\geq 1}$ named after him. 
These operators form a commutative 
algebra of operators. 
Hecke's student Petersson introduced a positive definite 
hermitian product on the space 
of cusp forms and the Hecke operators are hermitian for this product.
We thus can find a basis of eigenvectors, called eigenforms which are
eigenvectors for all the Hecke operators. Then one shows 
that for such a nonzero eigenform $f =\sum a(n) q^n$ 
with eigenvalue $\lambda(m)$ under $T(m)$ we have
$$
\lambda(m)a(1)=a(m)
$$
so that $a(1)\neq 0$
and when we divide $f$  by $a(1)$ to get $a(1)=1$ we see
that then the  Fourier coefficients are the eigenvalues. 

It was also Hecke who noticed that the fact that a modular form of 
weight~$k$ is an eigenform for
the Hecke algebra implies that the formal $L$-series of $f$ has an Euler
product:
$$
L(f,s)=\sum_{n\geq 1} \frac{a(n)}{n^s} =\prod_p \frac{1}{1-a(p)p^{-s}+p^{k-1-2s}}\, .
$$
For example, for the Eisenstein series one finds $L(E_{k},s)=\zeta(s)\zeta(s+1-k)$
with $\zeta(s)$ the Riemann zeta function.
Hecke also showed that the $L$-series of a cusp form can be extended
analytically to a holomorphic function of $s$ for all $s \in {\CC}$ and satisfies a 
functional equation, namely 
$$
\Lambda(f,s)=(-1)^{k/2}\Lambda(f,k-s)
$$ 
with $\Lambda(f,s)=\Gamma(s)(2\pi)^{-s} L(f,s)$.

Let $f$ be a (normalized with $a(1)=1$) cusp form of weight $k$ 
that is an eigenform of the Hecke algebra. We denote the totally real 
number field obtained by adjoining the eigenvalues of $f$ to ${\QQ}$ by
$K_f$. Deligne showed \cite{Deligne} that to such a cusp form 
one can associate an $\ell$-adic representation 
$$
\rho_{\ell,f}: {\rm Gal}(\overline{\QQ}/{\QQ})\to {\rm GL}(2,K_f\otimes{\QQ}_{\ell})
$$ of the Galois group of ${\QQ}$ 
with the property that for a prime $p$ different from $\ell$
the image $\rho_{\ell,f}(F_p)$ of a  geometric Frobenius element $F_p$ has 
trace $\lambda(p)$, the eigenvalue of $f$ under the operator $T(p)$,
 and determinant $p^{k-1}$. This was conjectured by Serre \cite{Serre-DPP}. 
This result exhibited the number theoretical meaning of modular forms.
The origin of these representations (namely the cohomology of modular varieties)
also enabled Deligne to apply his (later) result on the eigenvalues of Frobenius
on the cohomology of algebraic varieties to these eigenvalues, see \cite{Deligne}.
We shall see in Section \ref{CohomInt} where this representation lives.

\bigskip
  
\section{Siegel Modular Forms}
The notion of modular form (on ${\rm SL}(2,{\ZZ})$) generalizes. Consider a 
symplectic lattice ${\ZZ}^{2g}$, say with basis $e_1,\ldots,e_g,f_1,\ldots,f_g$,
which pair by $\langle e_i,e_j\rangle=0=\langle f_i,f_j\rangle$ and
$\langle e_i,f_j\rangle=\delta_{ij}$, the Kronecker delta. Then we let 
$$
\Gamma_g= {\rm Sp}(2g,{\ZZ})=
{\rm Aut}({\ZZ}^{2g},\langle \, , \, \rangle )
$$ 
be the automorphism group of this lattice, 
the so-called Siegel modular group. For $g=1$ one finds back 
${\rm SL}(2,{\ZZ})$. The action of ${\rm SL}(2,{\ZZ})$ on the upper half plane
generalizes to an action of  $\Gamma_g$  on the Siegel upper half space
$$
\mathfrak{H}_g=
\{ \tau \in {\rm Mat}(g \times g, {\CC}): \tau^t=\tau, {\rm Im}(\tau) >0\}
$$
via
$$
\tau \mapsto (a\tau+b)(c\tau+d)^{-1}
$$
for an element 
$$
\gamma = \left( \begin{matrix} a & b \\ c & d \\ \end{matrix} \right)
\in \Gamma_g
$$
where the decomposition of the $2g \times 2g$ matrix $\gamma$ in four $g\times g$ blocks corresponds to
the decomposition of ${\ZZ}^{2g}$ as a direct sum of two isotropic sublattices 
generated by the $e_i$ and the $f_i$. Often we shall write $\gamma=(a,b;c,d)$ for
such a matrix.

We can consider for $g>1$ 
holomorphic functions $f\colon  \mathfrak{H}_g \to {\CC}$ satisfying
$$
f((a\tau+b)(c\tau+d)^{-1})=\det(c\tau+d)^k f(\tau)
$$
for all $\tau \in \mathfrak{H}_g$ and 
$\gamma \in \Gamma_g$. This is the notion of a 
Siegel modular form of degree $g$ (or genus $g$) and weight $k$. 
But there is a wider generalization. Let
$$
\rho: {\rm GL}(g,{\CC}) \to {\rm Aut}(W)
$$
be a finite dimensional irreducible complex representation of ${\rm GL}(g,{\CC})$
on the vector space $W$. Then we can consider holomorphic maps
$$
f : \mathfrak{H}_g \to W
$$
satisfying
$$
f((a\tau+b)(c\tau+d)^{-1})=\rho(c\tau+d) f(\tau)
$$
for all $\tau \in \mathfrak{H}_g$ and $(a,b;c,d) \in \Gamma_g$. Such a function 
$f$ is called a Siegel modular form of degree $g$ and weight $\rho$. 
The space of all such modular forms is denoted by $M_{\rho}(\Gamma_g)$. 
If $\rho=\det{}^k$ 
we retrieve the notion of scalar-valued Siegel modular form of weight $k$.
It is a reflection of the Hartogs extension theorem that for $g>1$ we need not
demand a holomorphicity condition like we did for $g=1$ (when we demanded $a(n)=0$
for $n<0$). 
In fact, we can write the Fourier series of a Siegel modular form $f$
as
$$
f(\tau)=\sum_{n} a(n) q^n
$$
where the index $n$ runs over all $g\times g$ symmetric matrices whose entries
are half-integral and with integral diagonal entries and
where we employ the shorthand
$$
q^n = e^{2\pi i {\rm Tr} (n \, \tau)} 
$$
and $a(n)$ is a vector in the vector space $W$; 
the Koecher principle says
that if $a(n)\neq 0$ then $n$ is positive semi-definite. 

There is a way to obtain a Siegel modular form of degree $g-1$ from 
a form $f$ of degree $g$ by taking a limit
$$
\Phi(f)(\tau')= \lim_{t \to \infty} f(\begin{matrix} \tau' & 0 \\ 0 & i\, t 
\\ \end{matrix})
\quad \text{\rm with $t \in {\RR}_{>0}$ and $\tau' \in \mathfrak{H}_{g-1}$ .}
$$
This operator is called the Siegel operator.
Forms for which $\Phi(f)$ is zero are called cusp forms.

Like for $g=1$ we can form a graded ${\CC}$-algebra 
of scalar-valued Siegel modular forms
on $\Gamma_g$
$$
M(\Gamma_g)=\oplus_k M_k(\Gamma_g)
$$
where $M_k(\Gamma_g)=M_{{\det}^k}(\Gamma_g)$. The structure of this algebra is
known for $g=2$ and $g=3$ only. For $g=2$ Igusa proved
that the subalgebra  of forms of even weight is a polynomial algebra given by
$$
M^{\rm ev}(\Gamma_2)= \oplus_{k \equiv 0 (\bmod 2)} M_k(\Gamma_2)=
{\CC}[E_4,E_6,\chi_{10},\chi_{12}]\, .
$$
Here $E_4$ and $E_6$ are Eisenstein series, now for $g=2$, 
defined by series
$$
E_k=\sum_{(c,d)} \det(c\tau+d)^{-k}
$$
with a sum over non-associated pairs of co-prime symmetric 
integral matrices (where non-associated is taken with respect 
to left multiplication by elements of ${\rm GL}(2,{\ZZ})$),
and $\chi_{10}$ and $\chi_{12}$ are cusp forms; up to a 
normalization they can be given by
$E_{10}-E_4E_6$ and $E_{12}-E_6^2$. 
The algebra of modular forms of all weights
is an extension of $M^{\rm ev}(\Gamma_2)$ generated by a form
$\chi_{35}$ of weight $35$ that satisfies a quadratic relation 
$\chi_{35}^2=P(E_4,E_6,\chi_{10},\chi_{12})$ expressing it as a polynomial
in the even weight modular forms. For $g=3$ the structure of $M(\Gamma_3)$
was determined by Tsuyumine \cite{Tsuyumine}; it has $34$ generators 
and thus becomes   
rather complicated to handle.

Vector-valued Siegel modular forms have attracted much less attention than the
scalar-valued ones, but as we shall see they are the natural generalization of 
the elliptic modular forms. Much less is known about vector-valued Siegel
modular forms than about scalar-valued forms.

The quotient space ${\rm SL}(2,{\ZZ})\backslash \mathfrak{H}_1$ is the moduli
space of complex elliptic curves; 
similarly, ${\rm Sp}(2g,{\ZZ})\backslash \mathfrak{H}_g$ is the moduli space of
principally polarized complex abelian varieties of dimension $g$:
$$
{\mathcal A}_g({\CC}):= {\rm Sp}(2g,{\ZZ})\backslash \mathfrak{H}_g \, .
$$
It can be viewed as the fiber over the infinite place of a 
moduli stack ${\mathcal A}_g$ over ${\ZZ}$ parametrizing principally 
polarized abelian varieties and this
comes equipped with a universal abelian variety
$$
\pi: {\mathcal X}_g \longrightarrow {\mathcal A}_g
$$
with fibre $X_{\tau}={\CC}^g/{\ZZ}^g+\tau {\ZZ}^g$ over the orbit of $\tau$ in
${\mathcal A}_g({\CC})$.
We should view these moduli spaces (or quotients) as stacks or orbifolds.

Modular forms can be interpreted as sections of a vector bundle: $\Gamma_g$ acts
on $\mathfrak{H}_g \times {\CC}^g$ by
$$
(\tau,z) \mapsto (\gamma \tau, (c\tau +d)z)
$$
where $\gamma=(a,b;c,d)$. This defines a vector bundle (in the orbifold sense). 
More intrinsically, we have
$$
{\EE}_g= \pi_* \Omega^1_{{\mathcal X}_g/{\mathcal A}_g},
$$
the Hodge bundle,
and its restriction to ${\mathcal A}_g({\CC})$ 
corresponds to the standard representation of ${\rm GL}(g,{\CC})$.
For every irreducible representation $\rho$ of ${\rm GL}(g,{\CC})$ we have a corresponding
bundle $U_{\rho}$ on ${\mathcal A}_g({\CC})$ produced by applying an 
appropriate Schur functor to ${\EE}$; or more concretely, it can be 
described as the
quotient of the action of $\Gamma_g$ on ${\mathfrak H}_g \times W$ 
with $W$ the representation space of $\rho$ by
$(\tau,w) \mapsto (\gamma \tau, \rho(c \tau+d) w)$.
 Then we can interpret modular forms as sections:
for $g\geq 2$ we have
$$
M_{\rho}(\Gamma_g)=
\hbox{\rm space of sections of $U_{\rho}$ on ${\mathcal A}_g({\CC})$} \, .
$$
Our moduli space can be compactified by a toroidal compactification over which
the Hodge bundle extends. Then the bundles $U_{\rho}$ extend too.

By the Koecher principle these sections extend over a toroidal compactification
and we find that the space of cusp forms $S_{\rho}(\Gamma_g)$ consists of the
sections that vanish at infinity, that is, on the divisor that is added to
${\mathcal A}_g({\CC})$ to compactify it.

Like for $g=1$ one can define Hecke operators. They are induced by correspondences
and give a commutative algebra of operators. In particular for every prime $p$
we have an operator $T(p)$, but we have also operators $T_i(p^2)$ for $i=1,\ldots,g$.
We refer to \cite{Freitag, Andrianov}. Also the Petersson product has its analogue and 
we thus can find bases of eigenforms for the action of the Hecke algebra on the spaces of cusp forms.

\section{Cohomological interpretation}\label{CohomInt}

As alluded to in Section \ref{modforms} modular forms admit a cohomological
interpretation and this is a powerful tool as we shall see.
Let $\pi: {\mathcal X}_1 \to {\mathcal A}_1$ be the universal elliptic curve.
Then we find a local system ${\VV}=R^1\pi_*{\QQ}$ on ${\mathcal A}_1({\CC})$
of ${\QQ}$-vector spaces
with fibre $H^1(E,{\QQ})$ over $[E]$. This is a local system of rank $2$ 
associated to the standard representation of $\Gamma_1$. 
By taking ${\rm Sym}^a({\VV})$, the $a$th 
symmetric power of ${\VV}$, we obtain 
a local system ${\VV}_a$ of rank $a+1$.

This local system enters in the cohomological interpretation of elliptic modular cusp
forms that Eichler and Shimura 
gave in the 1950s.

\begin{theorem} (Eichler-Shimura)
For even $a\geq 2$ we have
$$
H^1_c({\mathcal A}_1 \otimes {\CC}, {\VV}_a\otimes {\CC})
\cong S_{a+2}(\Gamma_1) \oplus \overline{S}_{a+2}(\Gamma_1) \oplus {\CC}
$$
where $H^1_c$ stands for compactly supported cohomology and the right hand 
side displays the mixed Hodge structure of the left hand side. 
The space $S_{a+2}(\Gamma_1)$ has Hodge weight $(a+1,0)$ and $\overline{S}_{a+2}$
is the complex conjugate of $S_{a+2}$.
\end{theorem}

The case $a=0$ corresponds to the compactly supported cohomology of the moduli space ${\mathcal A}_1({\CC})$.
In this case $S_{2}=(0)$; but if we replace ${\rm SL}(2,{\ZZ})$ by a congruence subgroup $\Gamma_0(N)$
then cusp forms of weight $2$ might show up and these would give differential forms via 
$f(\tau) \mapsto f(\tau)d\tau$ and these would contribute to the cohomology of the quotient space
$\Gamma_0(N)\backslash \mathfrak{H}$ and its compactification
$\overline{\Gamma_0(N)\backslash \mathfrak{H}}$. 
This is how Eichler and Shimura were led to their result. 
Then the Galois representation associated to $f$ can be found
on the Tate module of the Jacobian of the compactification $\overline{\Gamma_0(N)\backslash \mathfrak{H}}$.

But the moduli space ${\mathcal A}_1$ is defined over ${\ZZ}$ and 
the Eichler-Shimura isomorphism has an $\ell$-adic counterpart, due to Deligne
\cite{Deligne}.
It is embodied in the identity
$$
{\rm Tr}(T(p), S_{a+2})=
{\rm Tr}(F_p, H^1_c({\mathcal A}_1\otimes \overline{\FF}_p, {\VV}_a^{(\ell)})) -1\, ,
\eqno(1)
$$
where ${\VV}^{(\ell)}$ denotes the $\ell$-adic sheaf $R^1\pi_* {\QQ}_{\ell}$
and ${\VV}_a^{(\ell)}$ its $a$th symmetric power,
and $F_p$ represents a geometric Frobenius element in the Galois group.
A convenient way to see where this identity comes from 
uses the Euler characteristic
$$
e_c({\mathcal A}_1, {\VV})= \sum (-1)^i [H_c^i({\mathcal A}, {\VV}_{a})]\, ,
$$
where we take the class $[H^i]$ of the cohomology group in an appropriate
Grothendieck group of mixed Hodge modules or of Galois representations
depending on whether we take complex cohomology or $\ell$-adic \'etale 
cohomology; then the statement is that
$$
e_c({\mathcal A}_1, {\VV}_{a})= -S[a+2]-1\, , \eqno(2)
$$
where $S[a+2]$ is the Chow motive defined by Scholl
associated to the space of cusp forms of weight $a+2$ on ${\rm SL}(2,{\ZZ})$.
The rank of this motive is $2s_{a+2}$ with $s_a=\dim S_{a+2}(\Gamma_1)$.
It is cut out by projectors on the cohomology of a symmetric power of the universal elliptic
curve over the moduli space and its compactification. 
Then the $\ell$-adic Galois representation associated to $S_{a+2}({\rm SL}(2,{\ZZ}))$
is to be found in 
$H^1_c({\mathcal A}_1 \otimes \overline{\FF}_p, {\VV}_a^{(\ell)})$; here we use
the comparison isomorphisms
$H^i_c({\mathcal A}_1\otimes \overline{\FF}_p, {\VV}^{(\ell)}_a) \cong
H^i_c({\mathcal A} \otimes \overline{\QQ}_p, {\VV}_a^{(\ell)})
\cong H^i_c({\mathcal A}\otimes \overline{\QQ}, {\VV}_a^{(\ell)})$ 
for $\ell\neq p$
and the surjection of Galois groups ${\rm Gal}(\overline{\QQ}_p/{\QQ}_p) \to
{\rm Gal}(\overline{\FF}_p/{\FF}_p)$.
\smallskip

The Hecke operators $T(n)$ 
are defined by correspondences on ${\mathcal A}_1$ and act on the cohomology.
Moreover, in characteristic $p$ there is the congruence relation that relates
the Hecke correspondence and the correspondence defined by Frobenius and its transpose.
The Hecke correspondence $T(p)$ is given by all pairs $(E,E')$ of elliptic 
curves linked by a cyclic isogeny $E \to E'$ of degree $p$. But in characteristic $p$
if $j$ and $j'$ are the $j$-invariants of such $E$ and $E'$ we have $j^p=j'$ or $j'^p=j$.

Then the identity (2) and the congruence relation 
imply that we can calculate the trace of the Hecke
operator $T(n)$ on $S_{a+2}$ from the trace of Frobenius on the cohomology
$H^1_c({\mathcal A}_1 \otimes \overline{\FF}_p, {\VV}_a^{(\ell)})$ via 
Deligne's identity (1).

\bigskip

From a calculational perspective 
this might seem pointless since we know generators for the
ring of modular forms and we can calculate the action of the Hecke operators
explicitly. Moreover, we have an explicit formula for the trace of the Hecke operator
$T(p)$, the Eichler-Selberg trace formula, see for example Zagier's appendix in
\cite{Lang}.
Nevertheless, on the contrary, our approach is a practical one. Indeed,
by the philosophy of Weil we can calculate cohomology by counting points over finite fields. 
And it turns out that we can use this to obtain information about 
modular forms. In practice we can do this by the following strategy:
\smallskip

\begin{enumerate}
\item{} make a list of all elliptic curves over ${\FF}_p$ up to isomorphism
over ${\FF}_p$;
\item{} for each elliptic curve $E$ in this list determine $\# E({\FF}_p)$
and $\# {\rm Aut}_{{\FF}_p}(E)$.
\end{enumerate}
We know by Hasse that the number of rational points of an elliptic curve 
over ${\FF}_p$ satisfies $\#E({\FF}_p)=p+1-\alpha-\overline{\alpha}$ with 
$\alpha$ an algebraic integer with $|\alpha|=\sqrt{p}$ and $\overline{\alpha}$
its complex conjugate.
Then we can calculate the trace of $T(p)$ on $S_{a+2}$ by
$$
1+ {\rm Tr}(T(p), S_{a+2})=
-\sum_{E} \frac{\alpha^a+\alpha^{a-1}\overline{\alpha}+\ldots + \overline{\alpha}^a}{\#{\rm Aut}_{{\FF}_p}(E)}\, ,
\eqno(3)
$$ 
where the sum is over the elliptic curves in our list. 
This formula is the concrete embodiment of Deligne's formula (1).
For example  for the Fourier coefficients of $\Delta$ we have
$$
\tau(p)=-1 - \sum_E \frac{\alpha^{10}+\alpha^9\overline{\alpha}+\ldots + \overline{\alpha}^{10}}{\#{\rm Aut}_{{\FF}_p}(E)}\, .
$$

Given a prime $p$ we can take the $E$ with the same trace $t$ together 
and then have to list how often a 
certain trace $t=\alpha+\overline{\alpha}$
occurs in our list, where we count  the frequency as
$$
w(t)=\sum_{E:\,  \#E({\FF}_p)=p+1-t} 
\frac{1}{\# {\rm Aut}_{{\FF}_p}(E)} \, ,
$$ 
where each $E$ with given trace counts with weight $1/\# {\rm Aut}_{{\FF}_p}(E)$.
For example for $p=17$ we thus get the list

\begin{footnotesize}

\smallskip
\vbox{
\bigskip\centerline{\def\quad{\hskip 0.6em\relax}
\def\quod{\hskip 0.5em\relax }
\vbox{\offinterlineskip
\hrule
\halign{&\vrule#&\strut\quod\hfil#\quad\cr
height2pt&\omit&&\omit&&\omit&&\omit&&\omit&&\omit&&\omit&&\omit&&\omit&&\omit&\cr
&$t$ && $\pm 8$ && $\pm 7$ && $\pm 6$ && $\pm 5$ && $\pm 4$ && $\pm 3$&& $\pm 2$ && $\pm 1 $&&$0$&\cr
\noalign{\hrule}
&$w$ && $1/4$ && $1/2$ && $3/2$ && $1/2$ && $1$ && $3/2$&& $7/4$&& $1/2$&&$2$& \cr
} \hrule}
}}
\end{footnotesize}

\smallskip

With this list we can calculate the trace of $T(17)$ 
on the space $S_k$ of cusp forms 
{\sl for all weights} $k\geq 4$ by formula (3).
Moreover, by having the list for a certain 
prime, computing the trace for an arbitrary 
weight is (almost) immediate. 
For prime powers $q$ we have a slightly modified formula.

As a final remark we note that the $1$ in formula (3) comes from the 
Eisenstein series $E_{a+2}$. Indeed, for even $a\geq 2$ we have
$$
H^1_{!}({\mathcal A}_1, {\VV}_{a})= S[a+2]
$$
where $H_{!}^i$ is the image of $H^i_c$ in $H^i$.
The $-1$ appearing in (1) and (2) 
comes from part of the eigenvalue $1+p^{a+1}$ of $E_{a+2}$ and it appears in
the kernel of the map $H^1_c \to H^1$; the other part $p^{a+1}$ appears in the
cokernel.

\bigskip
\section{Degree Two}
Since our knowledge about Siegel modular forms of degree $g\geq 2$ is much more limited
there is every reason to try to generalize this approach to higher~$g$. We consider then
$$
\pi: {\mathcal X}_g \to {\mathcal A}_g
$$
the universal abelian variety over our moduli space. This yields a local system
${\VV}=R^1\pi_*{\QQ}$ of rank $2g$ on ${\mathcal A}_g({\CC})$ 
and its $\ell$-adic variant ${\VV}^{(\ell)}$
on ${\mathcal A}_g\otimes {\CC}$ and on ${\mathcal A}_g \otimes {\FF}_{p}$.
The fibre over the point $[X]$ of a principally polarized abelian variety
$X$ is $H^1(X,{\QQ})$ or $H^1_{\rm et}(X,{\QQ}_{\ell})$. 
We shall simply write ${\VV}$ for ${\VV}^{(\ell)}$. From this local system we
can construct other local systems as follows. For every irreducible representation
of ${\rm Sp}(2g,{\QQ})$ of highest weight 
$\lambda=(\lambda_1\geq \lambda_2 \geq \cdots \geq \lambda_g)$ 
we get an associated local system ${\VV}_{\lambda}$; so with our conventions we
have
$$
{\VV}_{1,0,\ldots,0}={\VV} \, .
$$
We consider then the cohomology of this local system. As before, a convenient way 
to deal with it
is by using the so-called motivic Euler characteristic
$$
e_c({\mathcal A}_g({\CC}), {\VV}_{\lambda}\otimes {\CC})=
\sum_i (-1)^i [ H_c^i({\mathcal A}_g({\CC}), {\VV}_{\lambda}\otimes {\CC})]\, ,
$$
where $H^i_c$ is the compactly supported cohomology and the brackets refer to
the class in a Grothendieck group of mixed Hodge structures. 
But we can also consider $\ell$-adic \'etale cohomology 
and in this case we consider $H^i_c$ as a Galois representation and then $e_c$
lives in a Grothendieck group of Galois representations.

It is a result of Faltings \cite{F-C,Faltings}, but see also \cite{MM}, 
that $H^i$ and $H^i_c$ have a mixed Hodge filtration;
if we define the {\sl interior} cohomology by
$$
H^i_{!}=\hbox{\rm image of $H^i_c$ in $H^i$} \, ,
$$
then $H^i_{!}$ has a pure Hodge structure. Moreover, it follows from \cite{Faltings}
that if $\lambda$ is {\sl regular} 
(that is, $\lambda_1 > \lambda_2 > \cdots > \lambda_g$) then if $H^i_{!} \neq 0$
we must have $i=g(g+1)/2$, the dimension of ${\mathcal A}_g$.

Let us specialize to $g=2$. Then $\lambda$ is given as a pair $\lambda=(a,b)$ of
integers with $a\geq b \geq 0$ 
and according to Faltings we have a Hodge filtration
$$
F^{a+b+3} \subset F^{a+2} \subset F^{b+1} 
\subset F^0=H^3_{!}({\mathcal A}_2({\CC}), {\VV}_{a,b}\otimes {\CC}) \, .
$$
A main point is now that we have an interpretation of the first step:
$$
F^{a+b+3} \cong S_{a-b,b+3}, 
$$
with  $S_{a-b,b+3}$ the space of cusp forms on $\Gamma_2$ of weight $(a-b,b+3)$;
that is, these modular forms are the sections on $\tilde{\mathcal A}_2$ of
$$
{\rm Sym}^{a-b}({\EE}) \otimes \det({\EE})^{b+3} \otimes {\mathcal O}(-D)\, ,
$$
where $D$ is the divisor added to ${\mathcal A}_2$ to obtain an appropriate
 toroidal compactification $\tilde{\mathcal A}_2$ of ${\mathcal A}_2$.

\bigskip

Based on ample numerical evidence Carel Faber and I formulated a conjecture \cite{FvdG}
in 2004 that has now been confirmed completely. For regular local systems
($a>b>0$) this was done by Weissauer \cite{Weissauer} in 2009 and the irregular cases 
by Petersen \cite{Petersen} in 2013.

Before we state it we remark that it follows
 from the action of $-1_X$ on an abelian surface 
$X$ that if $a \not\equiv b \, (\bmod \, 2)$ then $e_c=0$. 
We thus may restrict ourselves to the case $a\equiv b \, (\bmod \, 2)$.

\begin{theorem}
We have
$$
{\rm Tr}(T(p),S_{a-b,b+3})=-{\rm Tr}(F_p,e_c({\mathcal A}\otimes \overline{\FF}_p,
{\VV}_{a,b})) +{\rm Tr}(F_p,e_{2,\rm extra}(a,b))
$$
with $e_{2,\rm extra}(a,b)$ a correction term given by
$$
s_{a-b+2}-s_{a+b+4}(S[a-b+2]+1)L^{b+1} + 
\begin{cases}
S[b+2] +1 & $a$ \quad even\\
-S[a+3] & $a$ \quad odd \\
\end{cases}
$$
\end{theorem}
Here $s_n=\dim S_n(\Gamma_1)$ and $L$ stands for the Lefschetz motive 
($=h^2({\PP}^1)$). And in order to make it work for $a=0$ we should define
$$ S[2]=-L-1 \quad \text{\rm and} \quad s_2=-1\, .$$

We arrived at this conjecture by applying the strategy mentioned 
for the case $g=1$:
first note that the moduli space of stable curves of genus $2$ of compact type
coincides with the moduli space ${\mathcal A}_2$ of principally polarized abelian
surfaces. We then make a list of all stable curves of compact type of genus 
$2$ over our finite field ${\FF}_q$ up to isomorphism over that finite field
and for each Jacobian $X$ of such a curve we compute the order of the group of
automorphisms defined over ${\FF}_q$ and the eigenvalues 
$\alpha_1, \bar{\alpha}_1,\alpha_2,\bar{\alpha}_2$ 
of Frobenius acting
on its $\ell$-adic Tate module 
(or on $H^1(X\otimes \overline{\FF}_q, {\QQ}_{\ell})$)
for a prime $\ell$ different from the characteristic. Then the trace of
Frobenius on  the cohomology of the local system ${\VV}_{a,b}^{(\ell)}$ is given by
summing over the $X$ in our list the values of a Schur function in the 
eigenvalues $\alpha_1, \bar{\alpha}_1,\alpha_2,\bar{\alpha}_2$,
that is, the $g=2$ analogue of (2).
For each $q$ this will give a number and we then tried to interpret this
number as a function of $q$ as a polynomial plus eigenvalue of an elliptic modular
form for small values of $a$ and $b$. This gave us an idea about the extra term
$e_{2,\rm extra}(a,b)$.

In order to make the list and in order to avoid the calculational problem
of determining the order of the automorphism group what we do is the following.
We take a suitable 
family of curves 
lying finitely over the moduli space and and then we count 
in this family and  divide the numbers by the degree of the map to the moduli space.

The way we arrived at the conjecture is a good example of experimental mathematics.
The fact that we had a dimensional formula for the spaces of Siegel modular forms
of degree $2$, due to Tsushima \cite{Tsushima}, was very helpful.
We also knew the numerical Euler characterics of the local systems.

If for given $q$ we have made such a list then this result enables us to compute
$
{\rm Tr}(T(q),S_{a-b,b+3}) 
$
for {\sl all values of} $(a,b)$ ! We did this for all $q<37$ (and have extended this now
to $q < 200$). 

Even for classical Siegel modular forms this gave information without reach before. 
For example, the space $S_{0,35}$ of scalar-valued Siegel modular forms of 
degree $2$ and weight $35$ is $1$-dimensional. It is generated by a cusp
form $\chi_{35}$ constructed by Igusa. It occurs in the cohomology of the 
local system ${\VV}_{32,32}$ where $e_{2,\rm extra}$  has the form
$$
e_{2,\rm extra}(32,32)=5\, L^{34}+S[34]
$$ 
and one finds as eigenvalue $\lambda(p)$ 
of $T(p)$ for $p\leq 37$:

\begin{footnotesize}
\smallskip
\vbox{
\bigskip\centerline{\def\quad{\hskip 0.6em\relax}
\def\quod{\hskip 0.5em\relax }
\vbox{\offinterlineskip
\hrule
\halign{&\vrule#&\strut\quod\hfil#\quad\cr
height2pt&\omit&&\omit&\cr
&$p$&& $\lambda(p)$  on $S_{0,35}$&\cr
height2pt&\omit&&\omit&\cr
\noalign{\hrule}
height2pt&\omit&&\omit&\cr
&$2$&&$-25073418240$&\cr
&$3$&&$-11824551571578840$&\cr
&$5$&&$9470081642319930937500$&\cr
&$7$&&$-10370198954152041951342796400$&\cr
&$11$&&$-8015071689632034858364818146947656$&\cr
&$13$&&$-20232136256107650938383898249808243380$&\cr
&$17$&&$118646313906984767985086867381297558266980$&\cr
&$19$&&$2995917272706383250746754589685425572441160$&\cr
&$23$&&$-1911372622140780013372223127008015060349898320$&\cr
&$29$&&$-2129327273873011547769345916418120573221438085460$&\cr
&$31$&&$-157348598498218445521620827876569519644874180822976$&\cr
&$37$&&$-47788585641545948035267859493926208327050656971703460$&\cr
height2pt&\omit&&\omit&\cr
} \hrule}
}}
\end{footnotesize}

Or consider the form $E_8 \chi_{35}$ that generates the $1$-dimensional space $S_{0,43}$. 
We have
$$
e_{2,{\rm extra}}(40,40)=S[42]+7\, L^{42}\, .
$$
Here we find the following eigenvalues $\lambda(p)$
of $T(p)$:

\begin{footnotesize}
\smallskip
\vbox{
\bigskip\centerline{\def\quad{\hskip 0.6em\relax}
\def\quod{\hskip 0.5em\relax }
\vbox{\offinterlineskip
\hrule
\halign{&\vrule#&\strut\quod\hfil#\quad\cr
height2pt&\omit&&\omit&\cr
&$p$&& $\lambda(p)$  on $S_{0,43}$&\cr
height2pt&\omit&&\omit&\cr
\noalign{\hrule}
height2pt&\omit&&\omit&\cr
&$2$&&$-4069732515840$&\cr
&$3$&&$-65782425978552959640$&\cr
&$5$&&$-44890110453445302863489062500$&\cr
&$7$&&$-19869584791339339681013202023932400$&\cr
&$11$&&$4257219659352273691494938669974303429235064$&\cr
&$13$&&$1189605571437888391664528208235356059600166220$&\cr
&$17$&&$-1392996132438667398495024262137449361275278473925020$&\cr
&$19$&&$-155890765104968381621459579332178224814423111191589240$&\cr
&$23$&&$-128837520803382146891405898440571781609554910722934311120$&\cr
&$29$&&$4716850092556381736632805755807948058560176106387507397101740$&\cr
&$31$&&$3518591320768311083473550005851115474157237215091087497259584$&\cr
&$37$&&$-80912457441638062043356244171113052936003605371913289553380964260$&\cr
height2pt&\omit&&\omit&\cr
} \hrule}
}}
\end{footnotesize}

\bigskip

For vector-valued modular forms it is just as powerful.
We illustrate this by giving the Hecke eigenvalues
for two vector-valued modular cusp forms, namely of weight $(14,7)$ and 
$(4,17)$. In both cases the dimension of the space of Siegel modular
cusp forms is $1$. We have
$$
e_{2,\rm extra}(18,4)=-L^5(S[16]+1)+2,
$$
and
$$
e_{2,\rm extra}(18,14)=-3L^{15}+S[16]+1\, .
$$
The graphs formed by the digits illustrate
the growth of the eigenvalues as roughly $p^{(j+2k-3)/2}$.
They also illustrate the fact that these eigenvalues tend to be `smooth'
(highly divisible) numbers for small values of $(j,k)$.

\begin{footnotesize}
\smallskip
\vbox{
\bigskip\centerline{\def\quad{\hskip 0.6em\relax}
\def\quod{\hskip 0.5em\relax }
\vbox{\offinterlineskip
\hrule
\halign{&\vrule#&\strut\quod\hfil#\quad\cr
height2pt&\omit&&\omit&&\omit&\cr
&$p$&& $\lambda(p)$  on $S_{14,7}$&&$\lambda(p)$ on $S_{4,17}$&\cr
height2pt&\omit&&\omit&&\omit&\cr
\noalign{\hrule}
height2pt&\omit&&\omit&&\omit&\cr
&$2$&&$-3696 $&&$-266112$&\cr
&$3$&&$511272$&&$-210323304$&\cr
&$5$&&$118996620$&&$668111687100$&\cr
&$7$&&$-82574511536$&&$-420920757352592$&\cr
&$11$&&$5064306707064$&&$-388201474991129976$&\cr
&$13$&&$-29379924792548$&&$28107151225966031596$&\cr
&$17$&&$170082580670244$&&$-3760611385645410867612$&\cr
&$19$&&$3752454431256520$&&$-6080023439267575397000$&\cr
&$23$&&$-79555863361862928$&&$168303583255503998515536$&\cr
&$29$&&$-81010055585118660$&&$ 15109310600861660971695180$&\cr
&$31$&&$-515521596253351616$&&$ 33344471645582702957462464$&\cr
&$37$&&$-40280723363343088436$&&$1247592679027009407366180988$&\cr
height2pt&\omit&&\omit&&\omit&\cr
} \hrule}
}}
\end{footnotesize}

Here we go only till $p\leq 37$, but in fact, 
we have extended the calculations to all $q< 200$.
The data will be made available on a website.

As an application of the availability of all these data we can mention 
the paper \cite{B-D-M} where all these data (traces for $q<150$) are used 
to approximate critical L-values used in 
checking congruences.

\bigskip

In the case of degree $2$ the analogue of the result of Deligne 
on the existence of a semi-simple $\ell$-adic 
Galois representation for a Hecke eigenform can be deduced 
from work of Taylor, Weissauer and Laumon
\cite{Taylor-Siegel3, Weissauer-book,Laumon}.
The dimension of the representation is $4$.
But unlike the situation for $g=1$ here it is not necessarily true that the four
eigenvalues of a Frobenius element at $p$ have a fixed absolute value $p^w$
for some $w$. When this was discovered this fact came as a surprise. 

In fact, at the end of the 1970s Kurokawa and Saito discovered (see \cite{Kurokawa-IM}) 
that there are
Siegel modular forms of degree $2$ that are lifts from elliptic modular forms.
In fact, if $f=\sum_n a(n) q^n$ is an elliptic cusp form of weight $2k-2$ with $k$ even
and a normalized (i.e.\ $a(1)=1$) 
eigenform for the Hecke algebra, then there is a scalar-valued 
Siegel modular form of even weight $k$ on $\Gamma_2$, also an eigenform for the 
Hecke algebra, with Hecke eigenvalues
$$
\lambda(p)= p^{k-2} +a(p)+p^{k-1} \, .
\eqno(4)
$$
Note that $a(p)=\alpha+\bar{\alpha}$ with $|\alpha|=p^{k-3/2}$.
So there is a $4$-dimensional Galois representation corresponding to 
(4) which we can realize as 
$$
{\QQ}(-k+2)\oplus R_f \oplus {\QQ}(-k+1)\, ,
$$
where $R_f$ is the $2$-dimension representation $\rho_{\ell,f}$ associated
to $f$ by Deligne.
But for reasons of cohomological weight only the $2$-dimensional part $R_f$ 
defined by $f$ can occur in the interior
cohomology of the local system ${\VV}_{k-3,k-3}$ on ${\mathcal A}_2$.
So the full $4$-dimensional Galois representation is not to be found in the cohomology
of the local system, only a part.

These lifts found by Kurokawa and Saito lie in a subspace of $S_k(\Gamma_2)$ called the Maass subspace
consisting of Siegel modular forms whose Fourier series
$$
\sum_{N \geq 0} a(N) e^{2\pi i {\rm Tr} (N \tau)}
$$
has the property that for $N=(n, r/2;r/2,m)$ the coefficient depends only
on the discriminant $d(N)=4mn-r^2$ and the content ${\rm g.c.d}(n,r,m)$. 
The $1-1$ correspondence between the Hecke eigenforms in $S_{2k-2}(\Gamma_1)$
and Hecke eigenforms in $S_{k}(\Gamma_2)$ was proved to a large extent 
by Maass, and completed by Andrianov and Zagier, see \cite{Z}.

\bigskip
In joint work with Bergstr\"om and Faber \cite{BFG}
we extended all this to the case of degree $2$ and level $2$. Here the group
$\Gamma_2$ is replaced by the kernel 
$\Gamma_2[2]=\ker( {\rm Sp}(4,{\ZZ}) \to {\rm Sp}(4,{\ZZ}/2{\ZZ}))$
of the natural reduction map. Since ${\rm Sp}(4,{\ZZ}/2{\ZZ})\cong \mathfrak{S}_6$,
the symmetric group on six letters, everything comes with an action of this group.
In \cite{BFG}  we formulated the conjectural 
analogue of the above theorem. The cohomology and spaces of modular forms
appear as representation spaces of $\mathfrak{S}_6$. Although it is
conjectural it is based on very firm evidence. 
In joint work with Fabien Cl\'ery and Sam Grushevsky \cite{C-vdG-G} 
we used these conjectural results
to determine the structure of a number of $M(\Gamma_2[2])$-modules
of vector-valued modular forms. 
By work of Igusa we know that the ring $M^{\rm ev}(\Gamma_2[2])$ 
of even-weight
scalar-valued  modular forms on $\Gamma_2[2]$ is generated by a $5$-dimensional
vector space of modular forms of weight $2$ that form a representation of type
$[2,1,1,1,1]$ for the group $\mathfrak{S}_6$ satisfying a quartic relation.
The ring of modular forms of all weights is a quadratic extension generated
by a modular form of weight $5$.
Then the direct sums
$$
\oplus_{k} M_{j,k}(\Gamma_2[2])
$$
for fixed value of $j$ are modules over $M(\Gamma_2[2])$. We determined 
the structure of such modules for small values of $j$.
The conjectural results suggested where to find generators of these modules
and by working with explicit generators and using the dimension formulas
the structure could be determined. 
We showed 
that the $M(\Gamma_2[2])$-module 
$$
\Sigma_2=\oplus_{k: \rm odd} S_{2,k}(\Gamma_2[2])
$$
is generated by ten modular forms $\Phi_i\in S_{2,5}$ ($i=1,\ldots,10$) 
with $\sum \Phi_i=0$ that span a $9$-dimensional
$\mathfrak{S}_6$-representation of type $[2,2,1,1]$. 
Since we know the decomposition
of $S_{2,k}$ as a $\mathfrak{S}_6$-representation we can read off the relations.
Similarly, the module
$$
\bigoplus_{k: \rm even} M_{2,k}(\Gamma_2[2])
$$
is generated by $15$ modular forms $G_{ij}\in S_{2,4}(\Gamma_2[2])$ 
with $1\leq i < j \leq 6$ that generate a $\mathfrak{S}_6$-representation
of type $[3,1,1,1]\oplus [2,1,1,1,1]$. 
\bigskip
 
Before we move to degree $3$ we remark that we
can formulate the result for genus $2$ differently, namely in the form
$$
e_c({\mathcal A}_2, {\VV}_{\lambda})= -S[a-b,b+3]+ e_{2,\rm extra}(a,b)
\eqno(5)
$$
in analogy with the case $g=1$, by assuming that there is a motive $S[a-b,b+3]$
of rank $4 \dim S_{a-b,b+3}$ with the property that
$$
\begin{aligned}
{\rm Tr}(T(p),S_{a-b,b+3})=& {\rm Tr}(F_p,S[a-b,b+3]) \\
=& -{\rm Tr}(F_p, e_c({\mathcal A}_2\otimes \overline{\FF}_p, {\VV}_{a,b})
+e_{2, \rm extra}(a,b) \\
\end{aligned}
$$
The motive $S[a-b,b+3]$ is still hypothetical, but would be the analogue of
the motive $S[a+2]$ constructed by Scholl. But even if we do not know whether
such a motive exists we can use $S[a-b,b+3]$ in the formulation above
just as a bookkeeping device with the property that 
$$
{\rm Tr}(T(p),S_{a-b,b+3})= {\rm Tr}(F_p,S[a-b,b+3])\, .
$$ 

\section{Degree Three}
Now we move to degree $g=3$. In this case the local systems are indexed by
$\lambda=(a,b,c)$ with integers $a\geq b \geq c\geq 0$. 
The goal is now to find an analogue 
of the formula (5); in other words the problem is to 
come up with a formula for $e_{3,\rm extra}(a,b,c)$.
In \cite{BFGg3}
we formulated a conjecture.

\begin{conjecture}
We have 
$$
e_c({\mathcal A}_3\otimes \overline{\FF}_p, {\VV}_{\lambda})=S[a-b,b-c,c+4]+e_{3, \rm extra }(a,b,c)
$$
with the correction term $e_{3, \rm extra }(a,b,c)$ defined by
$$
\begin{aligned}
e_{3, \rm extra}(a,b,c)= 
-e_c({\mathcal A}_2 \otimes \overline{\FF}_p,{\VV}_{a+1,b+1})
-e_{2, \rm extra}(a+1,b+1) \otimes S[c+2] & \\
+e_c({\mathcal A}_2 \otimes \overline{\FF}_p,{\VV}_{a+1,c})
+e_{2, \rm extra}(a+1,c) \otimes S[b+3] & \\
-e_c({\mathcal A}_2 \otimes \overline{\FF}_p,{\VV}_{b,c})
-e_{2, \rm extra}(b,c) \otimes S[a+4] & \\
\end{aligned}
$$
and where $S[a-b,b-c,c+4]$ denotes a hypothetical motive or bookkeeping device
such that
$$
{\rm Tr}(T(p),S_{a-b,b-c,c+4})={\rm Tr}(F_p,S[a-b,b-c,c+4]) \, .
$$
\end{conjecture}
Again we arrived at this conjecture by using the fact that in genus $3$ 
the moduli of principally polarized abelian varieties are close to the moduli of genus three curves. 
In fact, the map ${\mathcal M}_3 \to {\mathcal A}_3$ is of degree $2$ as a map of
stacks. Then we can calculate the trace of a power of Frobenius on the cohomology of
${\VV}_{\lambda}$ by making a list of all principally polarized abelian 
varieties of dimension $3$ over ${\FF}_q$ up to isomorphism over ${\FF}_q$
that are Jacobians of curves 
and calculating the order of the automorphism group and the eigenvalues of Frobenius. By summing a Schur function in these eigenvalues over this list (weighted by
the inverse orders of the automorphism groups) we find the trace. Again we 
then tried to interpret this number as a function of $q$ in terms of `known'
functions, like polynomials in $q$, traces of Hecke operators on spaces of
cusp forms of degree $1$ and $2$. That is how we arrived at the expressions 
for $e_{3,\rm extra}(a,b,c)$. In fact, 
the form of part of the Eisenstein cohomology
found in \cite{vdG2} guided our guesses and 
suggested a first approximation. Surprisingly the formula for degree $3$
looks simpler than the formula for degree $2$. We did not know the dimension of the
spaces of cusp forms, but we knew the numerical Euler characteristics $E_c$
(obtained by replacing the
terms in the conjecture by their dimensions) 
of the local systems by \cite{BvdG}. 

\smallskip

The evidence for this conjecture is considerable:
\begin{enumerate}
\item{} While we are summing rational numbers, since we divide by the order 
of the automorphism groups, the procedure 
produces integral values for the traces of the Hecke operators;
\item{} What we know about the dimensions of the spaces of modular forms agrees with the conjecture; 
for example,
if we know that the space of cusp forms vanishes we always find that
${\rm Tr}(e_c)={\rm Tr}(e_{\rm extra})$;
\item{} We know the numerical Euler characteristics $E_c$ by \cite{BvdG};  as it turns out
$E_c-E_{\rm extra}$ is always divisible by $8$.
\item{} In all vector-valued cases when $E_c=E_{\rm extra}$ (i.e.\ when 
we conjecture that the space of cusp forms vanishes),
we have ${\rm Tr}(e_c)={\rm Tr}(e_{\rm extra})$.
\item{} Recently Ta\"{\i}bi \cite{Taibi} has been able to calculate 
dimensions of Siegel modular forms
assuming the validity of Arthur's trace formula. His values fit ours.
\item{} There is consistency with the results of Chenevier and Renard \cite{C-R}
on the number of automorphic level one polarized algebraic regular automorphic 
representation of ${\rm GL}_n$ over ${\QQ}$ for $n\leq 8$.
\item{} Harder type congruences between elliptic modular forms 
and Siegel modular forms, 
see below.
\item{} M\'egarban\'e  
has computed a number of eigenvalues (up to $q=13$) 
for several Siegel modular forms of degree $3$ as they appear among 
the automorphic forms of ${\rm SO}(8)$. A number of consistency checks have been made. 
\item{} There is agreement with the eigenvalues computed in \cite{C-vdG2}.
\end{enumerate}

Let us give two examples. We have
$$
e_c({\mathcal A}_3,{\VV}_{11,5,2})=S[6,3,6] -S[12]L^3+L^7-L^3+1\, ,
$$
where $S[6,3,6]$ stands as above for a bookkeeping device. In this case
we should have 
$\dim S_{6,3,6}=1$.
The other example is for weight $(4,2,8)$. Again we should have 
$\dim S_{4,2,8}=1$. We have 
$$
e_c(10,6,4)=L^8+S[4,10]+S[4,2,8]-1 \, .
$$
The corresponding modular form is of type $G_2$. 
We have the following table of eigenvalues (where for powers of a prime
we use the convention of \cite{BFGg3}):

\smallskip
\vbox{
\bigskip\centerline{\def\quad{\hskip 0.6em\relax}
\def\quod{\hskip 0.5em\relax }
\vbox{\offinterlineskip
\hrule
\halign{&\vrule#&\strut\quod\hfil#\quad\cr
height2pt&\omit&&\omit&&\omit&&\omit \cr
&$p$&& $\lambda(q)$  on $S_{6,3,6}$&&$\lambda(q)$ on  $S_{4,2,8}$& \cr
height2pt&\omit&&\omit&&\omit&&\omit \cr
\noalign{\hrule}
height2pt&\omit&&\omit&&\omit&\cr
&$2$ && $0$&&$9504$& \cr
&$3$ && $-453600$&&$970272$& \cr
&$4$ && $10649600$&&$89719808$& \cr
&$5$ && $-119410200$&&$-106051896$& \cr
&$7$ && $12572892800$&&$112911962240$& \cr
&$8$ && $0$&&$1156260593664$& \cr
&$9$ && $-29108532600$&&$5756589166536$& \cr
&$11$&&$ -57063064032$&&$44411629220640$& \cr
&$13$ && $-25198577349400$&&$209295820896008$& \cr
&$16$ && $341411782197248$&&$-369164249202688$& \cr
&$17$ && $-107529004510200$&&$1230942201878664$& \cr
&$19$ && $1091588958605600$&&$51084504993278240$& \cr
height2pt&\omit&&\omit&&\omit&\cr
} \hrule}
}}

\bigskip
Our results for $g=3$ also tell what kind of lifts one finds.
We even could identify lifts that come from the Lie group $G_2$
as suggested by Gross and Savin. We refer to \cite[9.1]{BFGg3}.

\section{Harder Type Congruences}
Congruences between modular forms have a considerable history. 
There are well-known congruences between cusp forms and Eisenstein series 
like the congruence between the Hecke eigenvalues
$$
\tau(p) \equiv p^{11}+1 \, (\bmod \, 691)
$$
for $\Delta$ and the Eisenstein series of weight $12$.
Such congruences
occur for primes dividing the numerator of $B_k/2k$ with $B_k$ 
the $k$th Bernoulli number, see \cite{Deligne-Serre}. 
Swinnerton-Dyer determined all 
congruences modulo a prime for modular forms on 
${\rm SL}(2,{\ZZ})$, see \cite{Swinnerton-Dyer, Serre-SB}.

Kurokawa found several congruences between Siegel modular forms in 1979
in \cite{Kurokawa}. In particular, he conjectured the following. Suppose that a 
prime $\ell$ divides the
critical value $2k-2$ of the symmetric square $L$-series of an eigenform $f$ 
of weight $k$ on ${\rm SL}(2,{\ZZ})$. Then 
there should exist a cusp form $F$ of degree $2$, 
an eigenform of the Hecke operators, whose 
eigenvalues for $T(p)$  are congruent to those of the 
weight $k$ Klingen Eisenstein series 
modulo $\ell$ for every $p$. 
Some of these congruences
were proven by Mizumoto in 1986, see \cite{Mizumoto} and reconsidered
by Katsurada and Mizumoto \cite{Katsurada-Mizumoto, Dummigan}. 

In general if there is a Hecke invariant splitting of the vector space
of Siegel modular forms, one might expect congruences between the factors of
this splitting. Harder suspected congruences between cusp forms and Saito-Kurokawa lifts.
He was motivated by his work on the Eisenstein cohomology of $\Gamma_2$.

The moduli space ${\mathcal A}_g$ admits a `minimal' compactification 
${\mathcal A}_g^*$ 
obtained by 
mapping the quotient $\Gamma_g \backslash \mathfrak{H}_g$ into projective space
by using a basis of the space $M_k(\Gamma_g)$ of scalar-valued modular forms
of sufficiently high weight $k$ and then taking the closure. In other words,
by taking ${\rm Proj}$ of the graded ring of scalar-valued Siegel modular forms.
Set-theoretically we have
$$
{\mathcal A}_g^*={\mathcal A}_g \sqcup {\mathcal A}_{g-1} \sqcup\cdots 
\sqcup {\mathcal A}_0 \, .
$$
It is usually called the Satake or Baily-Borel compactification. 
This decomposition is mirrored in the
Siegel operator that associates to a Siegel modular form of degree $g$ a 
form of degree $g-1$ by the limiting procedure
$$
\Phi_g f(\tau') = \lim_{t \to \infty} f\left( \begin{matrix}   \tau' & 0 \\ 0 & it \\
\end{matrix} \right) \qquad \hbox{\rm with $\tau' \in \mathfrak{H}_{g-1}, t \in {\RR}_{> 0}$ .}
$$  
If $f\in M_k(\Gamma_g)$ then $\Phi_g(f) \in M_{k}(\Gamma_{g-1})$; if $f$ is a vector-valued
Siegel modular form of weight $\rho$ then $\Phi_g f$ is a form of weight $\rho'$
with $\rho'$ the irreducible representation of highest weight 
$(\lambda_1,\lambda_2,\ldots, \lambda_{g-1})$
if $\rho$ has highest weight $(\lambda_1,\lambda_2,\ldots, \lambda_{g})$.
The cusp forms are by definition the forms with $\Phi_gf=0$.
So the boundary of our moduli space is the place where modular forms
of different degrees interact. Using the interpretation with cohomology
one sees the interaction in the long exact sequences connecting cohomology on the
interior and on the boundary. Harder studied this in an extensive way using the
Borel-Serre compactification (a real manifold with corners) 
and Betti cohomology. He considers the cohomology
in the local system of ${\ZZ}$-modules ${\VV}_{\lambda}^{\ZZ}$ and over
a ring $R={\ZZ}[1/N]$ with an appropriate integer $N$ inverted, one
considers the long exact cohomology sequence
$$
H^3_c(\Gamma_2\backslash \mathfrak{H}_2,{\VV}_{\lambda}^{\ZZ} \otimes R) 
\to
H^3(\Gamma_2\backslash \mathfrak{H}_2,{\VV}_{\lambda}^{\ZZ} \otimes R) 
\to H^3(\partial(\overline{\Gamma_2\backslash \mathfrak{H}_2}),{\VV}_{\lambda}^{\ZZ} 
\otimes R) \, .
$$
Using the eigenvalues of the Hecke operators that have different absolute
values on $H^3_{!}$ and on the Eisenstein part, 
here taken with ${\QQ}$-coefficients, 
one sees that there is a direct sum decomposition of
$
H^3(\Gamma_2\backslash \mathfrak{H}_2,{\VV}_{\lambda})
$
as $H^3_{!}(\Gamma_2\backslash \mathfrak{H}_2,{\VV}_{\lambda})$
plus its image in 
$H^3(\partial(\overline{\Gamma_2\backslash \mathfrak{H}_2}),{\VV}_{\lambda} )$. This latter part contributes to the Eisenstein
cohomology and the classes there can be described by certain Eisenstein
series. For an Eisenstein series generating an  eigen space 
of the cohomology the constant term  is a 
product of critical values of $L$-functions associated to a modular form
on ${\rm SL}(2,{\ZZ})$. 
But if one works with integral cohomology (Harder uses Betti cohomology)
and intersects this with the direct sum of the interior
integral cohomology and the integral Eisenstein cohomology it might be
that this direct sum has a non-trivial index in the integral cohomology
$H^3(\Gamma_2\backslash \mathfrak{H}_2, {\VV}_{\lambda}^{\ZZ})$. This index
is related to primes dividing this critical values of the $L$ functions
occuring in the denominator of the constant term of this Eisenstein series.
This will lead to a congruence between a part of the interior cohomology and
part of the Eisenstein cohomology. In concrete terms it leads to a congruence
between the eigenvalues of an honest Siegel modular form of genus $2$ 
(that is, not a Saito-Kurokawa lift) and the eigenvalues of an elliptic
modular form. 

Harder had this idea already quite early, 
see his Lecture notes \cite{HarderLNM}, but it seemed difficult to check it.
We quote 
``Ich halte es f\"ur sehr interessant, numerische Rechnungen 
durch\-zu\-f\"uhren. .... Ich glaube, da{\ss}  die Implementierung 
hiervon auf einem Computer, einige schwierige Aufgabe darstellt. 
Sie wird viel Zeit und Rechenaufwand kosten..." (\cite[p.\ 101]{HarderLNM}).
But more than ten years later when it became feasible to calculate
the eigenvalues of examples of Siegel modular forms, Harder set himself the 
task of making such conjectures explicit. The first case he arrived at was
the conjectured congruence for an eigenform in $F$ in $S_{4,10}(\Gamma_2)$ and the
normalized elliptic cusp form $f=\sum a(n) q^n$ 
generating $S_{22}(\Gamma_1)$: if $\lambda(p)$ denotes the eigenvalue of $F$
for the Hecke operator $T(p)$ we should have
$$
\lambda(p)\equiv p^8 + a(p)+p^{13} \, (\bmod \, 41)  \eqno(6)
$$
for all primes $p$.
The eigenvalues we could calculate fitted perfectly, see Harder's report
on this in \cite{Ha}. Harder formulated a rather precise
conjecture and in many examples the congruences he predicted
could be confirmed numerically, see \cite{vdG1}. 
For other numerical checks see \cite{GRS}.
Some of the congruences have now been proved, 
for example the original congruence (6) modulo $41$ 
for weight $(4,10)$ by Chenevier and Lannes \cite{Chenevier-Lannes}
by their work on unimodular lattices of rank $24$ and the modular forms 
related to this. 

Harder's conjectures admit extensions to higher genus.
The underlying idea is that direct sum decompositions of rational cohomology
that are stable under the action of the Hecke algebra
do not necessarily hold for integral cohomology and this leads to
congruences. These congruences are predicted by divisibility
properties of critical values of $L$ functions associated to modular forms.

Since we have various Hecke invariant subspaces of the cohomology,
like the interior cohomology, Eisenstein cohomology, spaces of lifted
forms, we can expect in higher degree many different types of congruences. 

In our paper \cite{BFGg3} we formulated a number of extensions to
degree three Siegel modular forms. Let us give just one example. 
The conjecture says in this case that if we have eigenforms 
$f \in S_{c+2}(\Gamma_1)$ and $g \in S_{a+b+6}(\Gamma_1)$ and an 
ordinary prime  $\ell$ (in the field of eigenvalues of $f$ and $g$) dividing the critical value
at $s=a+c+5$ of the $L$-series belonging to ${\rm Sym}^2(f)\otimes g$,
then there should be a genuine Siegel eigenform 
$F \in S_{a-b,b-c,c+4}(\Gamma_3)$ such that 
$$
\lambda_F(p) \equiv \lambda_f(p) (p^{b+2}+\lambda_g(p)+p^{a+3}) \, \bmod \ell
$$
for all primes $p$. An example is the case $(a,b,c)=(13,11,10)$ with $F\in
S_{2,1,14}(\Gamma_3)$. The Euler characteristic is
$$
2 (L^{13} S[12]+1)-2L^{11}+S[2,1,14]
$$
and $\dim S_{2,1,14}=1$. 
The forms are $f=\Delta=q-24 q^2+ \ldots $ in $S_{12}$ and $g\in S_{30}$ 
with Fourier series $g=1+ (4320+96\sqrt{51349})\, q+ \ldots$ and
the prime dividing critical value of the $L$ series turns out to be
$\ell=199$.
The heuristic traces of $T(p)$ for $F\in S_{2,1,14}$ are in the table
below.
One checks that the congruence above is indeed satisfied for all $q\leq 17$,
for example for $p=2$ the norm of
$$
-2073600+24\, (2^{13}+4320+96\sqrt{51349}+2^{16})
$$
equals $-232402452480$ and
is divisible by $\ell=199$. 

\smallskip
\vbox{
\bigskip\centerline{\def\quad{\hskip 0.6em\relax}
\def\quod{\hskip 0.5em\relax }
\vbox{\offinterlineskip
\hrule
\halign{&\vrule#&\strut\quod\hfil#\quad\cr
height2pt&\omit&&\omit&&\omit&\cr
&$p$&& $\lambda(q)$  on $S_{2,1,14}$&\cr
height2pt&\omit&&\omit&\cr
\noalign{\hrule}
height2pt&\omit&&\omit&\cr
&$2$ && $-2073600$& \cr
&$3$ && $-1885952160$& \cr
&$4$ && $1080940298240$& \cr
&$5$ && $-26851408810200$& \cr
&$7$ && $-34909007533294720$& \cr
&$8$ && $-5890898142638899200$& \cr
&$9$ && $2339767572242234760$& \cr
&$11$ && $-660044916805998490272$& \cr
&$13$ && $-10848446812874015943640$& \cr
&$16$ && $461697465916451767451648$& \cr
&$17$ && $-2009932545573210270768120$& \cr
&$19$ && $17632053727783741943750240$& \cr
height2pt&\omit&&\omit&\cr
} \hrule}
}}

For further discussion of congruences between modular forms we refer to
the paper by Bergstr\"om and Dummigan \cite{B-D} where a vast generalization 
of Harder's conjecture is made. For finding this generalization the data
provided by our calculations were of great help.
\bigskip

\section{Picard Modular Forms}
The approach that we are advocating in this paper can be applied in other cases as well.
Shimura gave in 1964 a list of rational ball quotients \cite{Shimura}; these ball quotients
turn out to be moduli spaces of covers of the projective line and therefore
are amenable to our approach. Bergstr\"om and I have pursued this for two cases
of $2$-dimensional ball quotients. These ball quotients are Picard modular surfaces.

If $F$ is an imaginary quadratic field with ring of integers $O_F$ and discriminant
$D$ we consider the vector space $V=F^3$ with hermitian form $h$
$$
z_1\bar{z_2}+z_2\bar{z}_1+z_3\bar{z}_3
$$
and we let 
$$
G=\{ g \in {\rm GL}(3,F): h(gz,gw)=\eta(g) h(z,w)\}
$$
with $\eta(g) \in {\QQ}$ the so-called multiplier. Then $G$ is an algebraic
group defined over ${\QQ}$. If $G^0=\ker \eta$ we set $\Gamma= G^0({\ZZ})$
and $\Gamma_1= G^0({\ZZ}) \cap \ker{\det}$. An element $g \in G({\RR})$ with $\eta(g)>0$
acts on 
$$
B=\{ l \in V_{\CC}=V\otimes_{\QQ} {\RR}: \hbox{\rm $l$ a line with $h_{|l} <0$} \}\, ,
$$
the set of lines on which $h$ is negative definite. It is a complex $2$-ball
$$
B=\{(u,v) \in {\CC}^2: v+\bar{v}+u\bar{u} <0\} \, ,
$$
where $u=z_3/z_2$ and $v=z_1/z_2$.
The space $\Gamma \backslash B$ is an orbifold. It is not compact but can be
compactified by adding finitely points, called the cusps. These are singular 
and can be resolved by elliptic curves. These quotients are moduli spaces
of $3$-dimensional abelian varieties with endomorphisms from $O_F$. 

Picard studied these 
quotients and the corresponding modular forms 
already at the end of the nineteenth 
century, see \cite{Picard1,Picard2}.  Shimura took up their study again 
in the 1960s and in the 1970s Shintani wrote an unpublished paper
about the (vector-valued) modular forms associated to these groups. 
Later Shiga, Holzapfel, Feustel, Finis and others studied these 
$2$-dimensional spaces and the associated
scalar-valued modular forms, 
see \cite{Holzapfel, Finis}.
 
We looked at two cases in Shimura's list of ball
quotients; they have $F={\QQ}(\sqrt{-1})$ or $F={\QQ}(\sqrt{-3})$.

Let us consider the case $F={\QQ}(\sqrt{-3})$ and 
$$
\Gamma[\sqrt{-3}]=\{ g \in \Gamma: g \equiv 1 (\bmod \sqrt{-3})\}
$$
and
$$
\Gamma_1[\sqrt{-3}]=\{ g \in \Gamma[\sqrt{-3}], \det(g)=1\} \, .
$$
Then in this case the moduli space 
$\Gamma_1[\sqrt{-3}]\backslash B$ is a moduli 
space of curves: Galois covers of genus $3$ and degree $3$ 
with a marking of the
ramification points. These curves can be written as
$$
y^3=f(x)
$$
with $f$ a polynomial of degree $4$ with non-vanishing discriminant. 

The surface
$\Gamma_1[\sqrt{-3}]\backslash B$ can be compactified by adding four cusps.
It comes equipped with two vector bundles: we have an exact sequence
$$
0 \to L \longrightarrow V_{\CC} \longrightarrow Q \to 0
$$
where $L$ denotes the tautological line bundle and $Q$ the tautological quotient
of rank $2$.
With $U=Q^{\vee}$ we have the isomorphisms
$$
U \otimes L \cong \Omega^1(\log D), \qquad
L^3 \cong \Omega^2(\log D)\, ,
$$
where $D$ is the divisor that resolves the four cusps. Correspondingly, we have
two factors of automorphy
$$
j_1((u,v),g)=g_{21}u+g_{22} v + g_{23} \, ,
$$
corresponding to $L$, and
$$
j_2((u,v),g)= {\rm Jac}((u,v),g)^{-t}\cdot j_1((u,v),g)^{-1}
$$
with ${\rm Jac}$ the Jacobian of the action, corresponding to $U$.
This gives us the notion of a Picard modular form. 
A scalar-valued modular form
of weight $k$ is a holomorphic section of $L^k$ and a scalar-valued cusp form
of weight $k$ is a section of $\Omega^2 \otimes L^{k-3}$.
The ring of scalar-valued Picard modular forms in this case is known by work
of Shiga, Holzapfel and Feustel. We have
$$
M(\Gamma[\sqrt{-3}])=\oplus_k M_k(\Gamma[\sqrt{-3}])=
{\CC}[\varphi_0, \varphi_1,\varphi_2] \, ,
$$
where the $\varphi_i$  are modular forms of weight $3$ and
$$
M(\Gamma_1[\sqrt{-3}])={\CC}[\varphi_0, \varphi_1,\varphi_2,\zeta]/(\zeta^3-P(\varphi_0,\varphi_1,\varphi_2))
$$
with $\zeta$ a cusp form of weight $6$. 
Finis calculated in 1998 (see \cite{Finis}) some eigenvalues
under the Hecke operators
of low-weight Picard modular forms for $M(\Gamma_1[\sqrt{-3}])$. 

We also have vector-valued modular forms; these correspond to sections of
$$
{\rm Sym}^j(U) \otimes L^{k-3}
$$
Although Shintani's unpublished paper treats vector-valued Picard modular forms
there were no explicit examples in the literature before we considered this case.
Since our quotient parametrizes curves our 
cohomological approach also works here. We have a universal family
$$
\pi: {\mathcal X} \longrightarrow {\mathcal M}
$$
of curves of genus $3$ with an action of a group of order $3$.
The local system ${\VV}=R^1\pi_* {\QQ}$ is of rank $6$ and it or its 
\'etale $\ell$-adic variant splits after base change to $F$ as
$$
{\WW} \oplus {\WW}^{\prime} \, ,
$$
where ${\WW}$ (resp.\ ${\WW}'$) corresponds to the $\rho$-eigenspace 
(resp.\ $\rho^2$-eigenspace) 
under the action of the group of order $3$ and ${\WW}$ is of rank $3$. Note that after base change
we have $G \sim {\rm GL}(3,F) \times {\GG}_m$, where the multiplicative group
${\GG}_m$ corresponds to the multiplier $\eta$. The local systems
that we can consider here correspond to representations of this group.

We consider local systems ${\WW}_{\lambda}$ with $\lambda=(a+i,i,-b+i)$
$$
{\rm Sym}^a{\WW} \otimes {\rm Sym}^b{\WW}^{\prime} \otimes \det{\WW}^i\, ,
$$
where $\det{\WW}$ is a local system the 6th power of which is a constant
local system. So we may choose $i$ modulo $6$ (or even modulo $3$). 
We look at the motivic Euler characteristic
$$
e_c({\mathcal M}, {\WW}_{\lambda})= 
\sum_i (-1)^i [H_c^i({\mathcal M}, {\WW}_{\lambda})] \, ,
$$
again taken in an appropriate Grothendieck group 
of mixed Hodge structures or 
Galois representations. 
Let us denote as above the interior cohomology, that is, the image of compactly
supported cohomology in the usual cohomology, by $H^i_{!}$.
We know by work of Schwermer-Ragunathan that $H^i_{!}\neq (0)$ implies $i=2$
for regular $\lambda$ and
we can show that we have a
Hodge filtration on $H^2_{!}$ of the form
$$
F^{a+b+2} \subset F^{a+1} \subset F^0 = H^2_{!}({\mathcal M}, 
{\WW}_{\lambda})\, .
$$
As before, a main point is here that the first step 
$F^{a+b+2}$ can be interpreted as a space of cusp forms
$$
F^{a+b+2} \cong S_{b,a+3}(\Gamma[\sqrt{-3}], {\det}^i)\, ,
$$
namely, cusp forms of weight $(b,a+3)$ and character ${\det}^i$. Such
modular forms can be described as 
functions $f: B \to {\rm Sym}^j({\CC}^2)$ satisfying
$$
f(\gamma(u,v))={\det}(\gamma)^i\,  (j_1(\gamma,u,v))^{a+3} {\rm Sym}^b(j_2(\gamma,u,v))f(u,v)
$$
for $\gamma \in \Gamma[\sqrt{-3}]$
and a vanishing condition at the cusps.

We can make counts over finite fields of cardinality $q\equiv 1 (\bmod 3)$, and try to find the
extra term $e_{\rm extra}(a,b)$. It is composed of two parts: the Eisenstein cohomology and
the endoscopic part.

The Eisenstein cohomology  is given as
$$
e_{\rm Eis}({\mathcal M}, {\WW}_{\lambda})= \sum_i (-1)^i [{\ker}(H^i_c \to H^i ({\mathcal M}, {\WW}_{\lambda})] \, .
$$
This has been determined by Harder \cite{HarderU21, HarderLNM}. 
As a Galois representation it is given as an explicit sum of
Hecke characters of $F$. So we subtract it. 
But then the endoscopy remains. Indeed, there is a 
plethora of endoscopic terms. 
In \cite{JBvdG2} we determined heuristically all those terms. Subtracting these
gives us the trace of the Hecke operator on the spaces of 
honest Picard modular forms.
In fact we have a completely explicit but complicated 
conjectural formula for $e_c({\mathcal M}, {\WW}_{\lambda})$.  
We refer to \cite{JBvdG2} and content ourselves
by giving some simple examples. We note that the symmetric group
${\mathfrak S}_4$ acts on our modular surface (permuting the four cusps)
and we describe the cohomology as a representation
of the symmetric group ${\mathfrak S}_4$.
We denote an irreducible representation of ${\mathfrak S}_4$ by
 the corresponding partition, thus writing $s[4],\ldots, s[1,1,1,1]$. 

The Eisenstein cohomology vanishes for $\lambda=[a+i,i-b+i]$ with 
 $a\equiv b \equiv 2 (\bmod 3)$ and $i=1,2$.
Take $\lambda=[6k,0,0]$. Then
$$
e_{!}({\mathcal M}, {\WW}_{\lambda})= k (s[4]+s[3,1]+s[2,2]) L^{6k+1,1} +
S[0,6k+3,{\det}^2]
$$
Here $L^{\alpha,\beta}$ stands for a $1$-dimensional Galois representation of ${\rm Gal}(\overline{F}/F)$
or a Hecke character; moreover
$S[0,6k+3,{\det}^2]$ stands for a hypothetical motive or just as 
a bookkeeping device
with the property that
$$
{\rm Tr}(F_p, S[0,6k+3,{\det}^2])={\rm Tr}(T(p),S_{0,6k+3},\Gamma[\sqrt{-3}],{\det}^2)\, .
$$
Using this we can calculated traces of Hecke operators $T(p)$ for $p\leq 43$ 
equivariantly 
and the results fit for scalar-valued modular forms 
with the calculations done by Finis. And in fact this has now been extended
up to traces for $q<1000$.

We add an example involving vector-valued Picard modular forms: we let
$$
\lambda=(6k+1,0,-1)
$$
and in this case we find the conjectural formula
$$
e_{!}({\mathcal M}, {\WW}_{\lambda})=(ks[4]+ks[3,1]+(k+1)s[2,2])L^{6k+2,2}
+
S[1,6k+4,{\det}^2]
$$ 
We are still in the process of writing up these results. But in the 
meantime Fabien Cl\'ery and I
have used these heuristic results 
to investigate  modules  of vector-valued Picard modular forms
over the ring of scalar-valued Picard modular forms \cite{C-vdG1}. 
For example, consider 
$$
\Sigma_1^{(i)}= \oplus_{k=1}^{\infty} S_{1,3k+1}(\Gamma[\sqrt{-3}],{\det}^i) \qquad \hbox{\rm for $i=0,1,2$}.
$$
This is a module over the ring of scalar-valued modular forms
$$
M=\oplus_k M_{0,3k}(\Gamma[\sqrt{-3}])={\CC}[\varphi_0,\varphi_1,\varphi_2]\, .
$$
We proved the following result.

\begin{theorem}
The module $\Sigma_1^{0}$ is generated over the ring $M$ of scalar-valued
Picard modular forms by three vector-valued modular forms
$\Phi_0,\Phi_1,\Phi_2 \in S_{1,7}(\Gamma[\sqrt{-3}])$ with the relation
$$
\varphi_0 \Phi_0 + \varphi_1 \Phi_1 +\varphi_2 \Phi_2 =0\, .
$$
\end{theorem}

The modular forms $\Phi_i$ appearing here 
are constructed as a sort of Cohen-Rankin 
bracket, with for example $\Phi_2$ given by
$$
\Phi_2= [\varphi_0,\varphi_1] \sim 
\varphi_0 \left( 
\begin{matrix} \partial\varphi_1/\partial u \\
\partial \varphi_1 /\partial v \\ 
\end{matrix} \right) 
- \varphi_1 \left( \begin{matrix} \partial\varphi_0/\partial u \\
\partial \varphi_0 /\partial v \\ \end{matrix} \right) \, .
$$

We can calculate the eigenvalues for these forms and the 
values fit the heuristic values perfectly. The description of the module 
makes it possible to construct explicitly eigenforms and calculate
eigenvalues and these have been used to confirm the conjectural
formulas for $e_c({\mathcal M}, {\WW}_{\lambda})$ of \cite{JBvdG2}.
We just give one table that illustrates the case of the forms $\Phi_i$.

\begin{footnotesize}
\smallskip
\vbox{
\bigskip\centerline{\def\quad{\hskip 0.6em\relax}
\def\quod{\hskip 0.5em\relax }
\vbox{\offinterlineskip
\hrule
\halign{&\vrule#&\strut\quod\hfil#\quad\cr
height2pt&\omit&&\omit&&\omit&\cr
&$\alpha$ && $p$ && $\lambda_{\alpha}(\Phi_i)$ &\cr
\noalign{\hrule}
&$1+3\rho$&& $7$ && $759+261\rho$ & \cr
&$1-3\rho$&& $13$ && $-4137+1683\rho$ & \cr
& $-2+3\rho$ && $19$ && $24042+14733\rho$ & \cr
& $1+6\rho$ && $31$ && $-145401-241830\rho$ & \cr
& $4-3\rho$ && $37$ && $12900-114849\rho$  & \cr
& $1-6\rho$ && $43$ && $246567-8946\rho$  & \cr
& $4+9\rho$ && $61$ && $1048836-173205\rho$ & \cr
& $-2-9\rho$ && $67$ && $-1539510-1246887\rho$ & \cr
& $1+9\rho$ && $73$ && $-1563729+1261143\rho$  & \cr
& $7-3\rho$ && $79$ && $9921297+3294171\rho$  & \cr
& $-8+3\rho$ && $97$ && $5678616-3870891\rho$  & \cr
& $-2$ && $2$ && $72$ & \cr
& $-5$ && $5$ && $89622$ & \cr
} \hrule}
}}
\end{footnotesize}

\bigskip

For the details and the structure of $\Sigma_1^{1}$ and
$\Sigma_1^{2}$
and various other modules we refer to \cite{C-vdG1}.

\bigskip


\end{document}